\documentclass{amsart}

\usepackage{amsmath}
\usepackage{amsfonts}
\usepackage{amssymb}
\usepackage{latexsym}
\usepackage{epsfig}
\usepackage{psfrag}

\def\co{\colon\thinspace}
\newcommand{\begriff}[1]{\textbf{#1}}
\newcommand{\homeo}{\approx}
\renewcommand{\d}{\partial}
\renewcommand{\t}{\mathcal T}

\newtheorem{theorem}{Theorem}

\newtheorem{lemma}{Lemma}

\theoremstyle{definition}
\newtheorem{definition}{Definition}

\theoremstyle{remark}

\begin{document}

%*******************************************************************

\title{Almost normal Heegaard surfaces}  
\author{Simon A. King}
\address{Department of Mathematics,
Darmstadt University of Technology,
Schlossgartenstr.~7,
64289 Darmstadt,
  Germany}
\email{king@mathematik.tu-darmstadt.de}
\begin{abstract}
  We present a new and shorter proof of Stocking's result that any
  strongly irreducible Heegaard surface of a closed orientable
  triangulated $3$--manifold is isotopic to an almost normal surface.
  We also re-prove a result of Jaco and Rubinstein on normal spheres.
  Both proofs are based on the ``reduction'' technique
  introduced by the author.
\end{abstract}

 %Niederdimensionale Themen, die sonst nirgends passen


\maketitle

%*******************************************************************

\section{Introduction}
\label{sec:intro}

%%%%%%%%%%%% Inhaltliche Intro %%%%%%%%%%%%%%%%%%%%%%

Any closed orientable triangulated $3$--manifold $M$ has a
\begriff{Heegaard surface}, which is an embedded surface that decomposes
$M$ into two handlebodies.
By adding trivial handles, one can construct Heegaard surfaces of $M$ of
arbitrarily large genus. An important problem is to construct a Heegaard
surface of $M$ of \emph{minimal} genus.  Though this problem has an
algorithmic solution for closed orientable atoroidal Haken manifolds by
work of Johannson~\cite{johannson}, it is still open for non-Haken
manifolds.
One approach to the non-Haken case is based on a result of Casson and
Gordon~\cite{cassongordon}. It states that any minimal genus Heegaard
surface $H$ of a closed orientable irreducible non-Haken manifold $M$
is \begriff{strongly irreducible}, i.e., any two simple closed
essential curves on $H$ bounding embedded discs in different connected
components of $M\setminus H$ intersect each other in at least two
points.
In fact {any} compact orientable $3$--manifold can be decomposed
along incompressible surfaces so that the pieces have strongly
irreducible Heegaard surfaces with boundary~\cite{scharlthompson}.

Stocking~\cite{stocking} has shown that any strongly irreducible Heegaard
surface is isotopic to a so-called almost normal surface, a mild
generalization of normal surfaces that goes back to  Rubinstein~\cite{rubinstein}.
This is amazing, since normal surfaces, introduced by
Kneser~\cite{kneser} in his study of connected sums of $3$--manifolds,
have been designed to deal with incompressible surfaces, whereas
Heegaard surfaces bound two handlebodies and are thus completely
compressible on both sides.
The aim of this paper is to present a new and shorter proof
of Stocking's result.

%%%%%%%%%%%%%%% Technische Intro %%%%%%%%%%%%%%%%%%%%

\begin{figure}[htbp]
  \begin{center}
    \leavevmode
    \epsfig{file=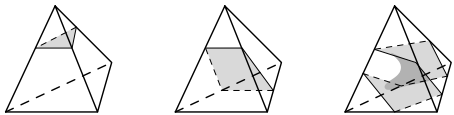,width=6cm}
    \caption{A triangle, a square and an octagon}
    \label{fig:pieces}
  \end{center}
\end{figure}
We follow in this paper the terminology of Matveev~\cite{matv}. We
define it in detail in Section~\ref{sec:almost} and here give just an
outline.
Let $M$ be a closed orientable $3$--manifold with a triangulation
$\t$.  A closed embedded surface in $M$ is $2$--normal if its
intersection with any tetrahedron is formed by copies of the pieces
shown in Figure~\ref{fig:pieces}, which are called normal triangles,
squares and octagons.
The surface is $1$--normal if it is formed by normal triangles and
squares only.

Let $\gamma$ be an \begriff{unknotted arc} for some $1$--normal
surface $F\subset M$, i.e., $\gamma$ is contained in a tetrahedron
$t$; connects two different components of $F\cap t$; and there is an
embedded disc in $t$ whose boundary consists of $\gamma$, one arc in
$\d t$ and two arcs in $F\cap t$. From $F$ we cut out small discs 
around the endpoints of $\gamma$, glue in an annulus (a ``tube'')
along $\gamma$ and denote the result by $F^\gamma$ --- see
Section~\ref{sec:basics} for a more general definition.

\begin{theorem}\label{thm:heegaard-an}
  Let $M$ be a closed orientable irreducible triangulated
  $3$--manifold, and let $H\subset M$ be a strongly irreducible
  Heegaard surface. Either there is a $1$--normal surface $F\subset M$
  and an unknotted arc $\gamma$ such that $F^\gamma\simeq H$, or $H$
  is isotopic to a $2$--normal surface with exactly one octagon.
\end{theorem}

We briefly outline our proof. It is easy to show that there are
$1$--normal surfaces $F_-,F_+$ that can be ``completed'' to form disjoint
surfaces $S_-,S_+\simeq H$, by adding some tubes disjoint from the edges
of $\t$ (possibly multiple, knotted, nested tubes). Here $S_-$ is a so-called
upper trivial and $S_+$ a lower trivial surface.
By application of a classical finiteness result of Kneser for
$1$--normal surfaces, we choose $F_-$ and $F_+$ so that any $1$--normal
surface with an upper (resp. lower) trivial completion isotopic to $H$
separating $S_-$ from $S_+$ is ``essentially the same'' as $F_-$
(resp. $F_+$).

If $S_-$ (or similarly $S_+$) satisfies a certain technical condition
then we apply the reduction technique that we have developed
in~\cite{king1} and~\cite{king2}. If the first alternative of
Theorem~\ref{thm:heegaard-an} does not hold then the result of the
reduction is isotopic to $H$ and is an upper trivial completion of
some $1$--normal surface different from $F_-$ (in contradiction to
the choice of $F_-,F_+$).

If, instead, $S_-$ and $S_+$ do not satisfy the mentioned
technical condition, we apply a \emph{thin position} argument. If
the first alternative of Theorem~\ref{thm:heegaard-an} does not hold, then
it yields a surface $S\simeq H$ that is an upper or lower trivial
completion of a $2$--normal surface with exactly one octagon.
Again, we apply the reduction technique to $S$, and derive a
contradiction to the choice of $F_-,F_+$, unless the second alternative of
Theorem~\ref{thm:heegaard-an} holds. Thus in either case, one
alternative is true.

The main reason why our proof is shorter than Stocking's is the
application of the completion and reduction techniques. This makes thin
position a more efficient tool. Actually we only need to know one
basic fact about thin position: A level surface in a thin position
embedding has no pair of nested or independent compressing discs. 
The reduction technique was used in~\cite{king1} to obtain bounds for
the bridge number of links formed by edges of a triangulation of
$S^3$. It also allows to simplify the correctness proof of the
Rubinstein--Thompson algorithm (\cite{thompson}, \cite{matv}), which
is implicit in~\cite{king1}.

We took a few technical results from the literature. To make this
paper essentially self-contained, we include proofs or at least proof
sketches for most cited results. Specifically, we re-prove
Scharlemann's ``no nesting'' lemma~\cite{scharlemann}, that is at the
base of our completion technique.
It is convenient for us, though not strictly necessary, to work in a
submanifold $N\subset M$ whose complement is a $3$--ball containing
all vertices of $\t$, and any $1$--normal sphere in $N$ is a copy of
$\d N$.  It is known~\cite{jacorubinstein} that such a submanifold
exists for any triangulated irreducible $3$--manifold different from
$S^3$.
We include here a proof that is another application of the reduction
technique. An algorithm that can be used to construct the submanifold
appears in~\cite{king1}, with a bound for the complexity of $\d N$.

The paper is organized as follows. In
Sections~\ref{sec:basics}--\ref{sec:uc-normal}, we expose some basic
constructions and the notions of almost $k$--normal and impermeable
surfaces.
In Section~\ref{sec:reductions}, we recall the reduction technique and
prove the existence of the above mentioned submanifold $N$.
In Section~\ref{sec:completion}, we introduce upper and lower
completions of surfaces and re-prove the ``no nesting'' lemma of
Scharlemann. Theorem~\ref{thm:heegaard-an} is proved in
Section~\ref{sec:proof}.

%*******************************************************************
%*******************************************************************

\section{Notations and basic notions}
\label{sec:basics}

In this paper, $M$ denotes a closed orientable $3$--manifold with a
triangulation $\t$.  We denote the $i$--skeleton of $\t$ by $\t^i$, for
$i=0,\dots,3$. An ambient isotopy that fixes each simplex of $\t^i$
setwise is an \begriff{isotopy mod $\t^i$}.
We denote the number of connected components of a topological space
$X$ by $\#(X)$. The notation $X\subset M$ stands for a tame embedding of
$X$ into $M$, and $U(X)$ denotes an open regular neighborhood of $X$
in $M$.

Let $S\subset M\setminus \t^0$ be a closed surface, which is allowed
to be empty or non-connected. The \begriff{weight} of $S$ is $\|S\|=\#(S\cap
\t^1)$.  The surface $S$ is \begriff{splitting}
if any connected component of $S$ decomposes $M$ into two parts.
We fix a vertex $x_0\in \t^0$ and define $B^+(S)$ (resp.\ $B^-(S)$) as
the closure of the union of components of $M\setminus S$ that are
connectable with $x_0$ by a path transversely intersecting $S$ in an
odd (resp.\ even) number of points.  In particular, $B^-(\emptyset)=M$.
We do not include $x_0$ in the notation ``$B^+(S)$'', since in our
applications the choice of $x_0$ plays no role.

Let $\Gamma\subset B^+(S)$ be an embedded graph, not necessarily
transversal to $\t^2$. We define
$S^\Gamma= \d(B^-(S)\cup U(\Gamma))$. Similarly, if $\Gamma\subset
B^-(S)$ then we define $S^\Gamma = \d(B^+(S)\cup U(\Gamma))$.

If $S$ is in general position to $\t$ then $G(S)$ denotes the union
of those connected components of $S\cap \t^2$ that intersect $\t^1$. For
any tetrahedron $t$ of $\t$, $G(S)\cap \d t$ is a disjoint union of
circles.
Hence by inserting disjoint discs in $M\setminus \t^2$, one obtains
from $G(S)$ a closed embedded surface $S^\times$ (pronounced ``$S$ cut''),
so that $S^\times\cap \t^2 =G(S)\subset S\cap \t^2$ and $S^\times
\setminus \t^2$ is a disjoint union of discs. One obtains $S^\times$
from $S$ by cut-and-paste operations along discs in $M\setminus \t^2$,
and omission of connected components that are disjoint from $\t^1$.

%*********************************************************************

\section{Almost $k$--normal surfaces}
\label{sec:almost}

An embedded arc $\gamma$ in a closed $2$--simplex $\sigma$ of $\t$
disjoint from the vertices of $\sigma$ with
$\gamma\cap\d\sigma=\d\gamma$ is a \begriff{normal arc}, if $\gamma$
connects different edges of $\sigma$.  Otherwise, $\gamma$ is a
\begriff{return}. 
%
%\begin{definition}\label{def:akn}
  A closed surface $S\subset M$ in general position to $\t^2$ is
  \begriff{almost $k$--normal}, if $S\cap \t^2$ is a union of normal
  arcs and of circles in $\t^2\setminus\t^1$, and for any tetrahedron
  $t$ of $\t$, any edge $e$ of $t$ and any connected component $c$ of
  $S\cap \d t$ holds $\#(c\cap e)\le k$.
%\end{definition}
%

Almost $k$--normal surfaces are difficult to deal with. For instance,
if an almost $k$--normal surface meets a tetrahedron in an annulus
then this annulus can be arbitrarily knotted.
These problems do not occur if one makes the additional assumption
that all intersections with tetrahedra are discs, as in the following
definition.
\begin{definition}
  Let $F\subset M$ be an almost $k$--normal surface. 
  If $F\setminus\t^2$ is a disjoint 
  union of discs and $F\cap\t^2$ is a union of normal arcs in the
  2--simplices of $\t$ then $F$   is a \begriff{$k$--normal surface}. 
\end{definition}

Note that we do not assume $F\not=\emptyset$. Actually it is
convenient in our proofs to consider the empty set as a $1$--normal
surface.
If $S\subset M$ is an almost $k$--normal surface, then $S^\times$ is
$k$--normal. It is easy to show that $S^\times$ is  determined
by $S\cap \t^1$ up to isotopy mod $\t^2$. We have $S^\times =\emptyset
$ if and only if $\|S\|=0$.

%***********************

In this paper, we mainly consider $1$-- and $2$--normal surfaces.  A
$2$--normal surface intersects a tetrahedron $t$ of $\t$ in a disjoint
union of {triangles}, squares and octagons, as in
Figure~\ref{fig:pieces}.  Up to isotopy mod $\t^2$ there are four
types of triangles (one for each vertex of $t$), three types of
squares (one for each pair of opposite edges of $t$), and three
types of octagons in $t$. A $1$--normal surface is formed
by triangles and squares.
\begin{theorem}[Lemma~4 in~\cite{haken2}]\label{thm:kneserhaken}
  Let $n$ be the number of tetrahedra of $\t$.
  Let $F\subset M$ be a $1$--normal surface with more than
  $20 n$ connected components. Then two connected components of $F$ 
  are isotopic mod $\t^2$.\qed 
\end{theorem}
This finiteness result goes back to Kneser~\cite{kneser}.  There are
better bounds for the number of components, but the proofs are more
complicated than the proof of Lemma~4 in~\cite{haken2}.

%*********************************************************************

\section{Compressing and essential discs} 
\label{sec:uc-normal}

Here we study some properties of embedded discs whose interior is disjoint
from $\t^1$. In the next section, these discs will be
used to construct isotopies of $S$.

\begin{definition}
  Let $S\subset M$ be a closed surface in general position to $\t^1$. Let
  $\alpha\subset\t^1\setminus\t^0$ and $\beta\subset S$ be embedded
  arcs with $\d\alpha=\d \beta$. A compact embedded disc $D\subset M$
  in general position to $S$ is a \begriff{compressing disc} for $S$ with
  \begriff{string} $\alpha$ and \begriff{base} $\beta$, if $\d
  D=\alpha\cup\beta$ and $D\cap \t^1=\alpha$.
\end{definition}

\begin{definition}\label{def:essential}
  Let $S\subset M$ be a closed surface in general position to $\t^1$. A compact
  embedded disc $D\subset M\setminus \t^2$ is
  \begriff{essential} for $S$, if
  $\d D\subset S$ is not null-homotopic in $S\setminus \t^2$ and
  $\#(D\cap S)\le \#(D'\cap S)$ for any compact embedded disc
  $D'\subset M\setminus \t^2$ bounded by $\d D$.
\end{definition}

Let $D$ be a compressing or essential disc for a splitting surface
$S\subset M$. If $\d D\cap S$ has a collar in $D\cap B^+(S)$ (resp.\ 
$D\cap B^-(S)$) then $D$ is an \begriff{upper} (resp.\ 
\begriff{lower}) compressing or essential disc.  If $D\cap S\subset \d
D$ then $D$ is \begriff{strict}.

Obviously, $k$--normal surfaces have no essential discs. In fact, as
stated in the following lemma, all discs in the complement of $\t^1$
with boundary in a $1$--normal surface are trivial.
\begin{lemma}[Lemma~10 and~11 in~\cite{king1}]\label{lem:discparallel}
  A $1$--normal surface $F\subset M$ has no compressing discs. Let
  $D\subset M\setminus \t^1$ be a closed embedded disc in general
  position to $F$ with $\d D\subset F$.  Then $\d D$ bounds a disc in
  $F\setminus\t^1$.
\end{lemma}
\begin{proof}
  We choose the disc $D$ up to isotopy of the pair $(D,\d D)$ in
  $(M\setminus \t^1, F\setminus \t^1)$ so that $\#(\d D\cap \t^2)$ is
  minimal. By cut-and-paste arguments, there are no circles in $D\cap
  \t$. If $\d D$ is not contained in a single tetrahedron then there
  is a component $D'$ of $D\setminus \t^2$ such that $\d D'\cap F$ is
  a single arc, contained in some tetrahedron $t$.
  Since $F$ is $1$--normal, $\d D'\cap F$ splits off a disc
  $D''\subset t$ from some component of $F\cap t$ so that $\d D''\cap
  \t^1=\emptyset$.  By an isotopy of $D$ with support in $U(D'')$, one
  can remove $\d D$ from $D''$, decreasing $\#(\d D\cap \t^2)$ --- in
  contradiction to our choice.
  Therefore $\d D$ is contained in a single tetrahedron, thus bounding
  a disc in $F\setminus \t^2$. The arguments are similar when $D$ is a
  compressing disc (see~\cite{king1} for details).
\end{proof}

In the remainder of this section, we recall the notion of impermeable
surfaces and its relationship with almost $2$--normal surfaces.
Let $S\subset M$ be a splitting surface in general position to $\t$
and let $D_1, D_2$ be upper and lower compressing discs for $S$ with
strings $\alpha_1,\alpha_2$. If $D_1\subset D_2$ or $D_2\subset D_1$
then $D_1$ and $D_2$ are \begriff{nested}. If $D_1\cap D_2\subset
\d\alpha_1\cap \d\alpha_2$ then $D_1$ and $D_2$ are
\begriff{independent} from each other.
If $S$ has both strict upper and strict lower compressing discs and has
no pair of nested or independent upper and lower compressing discs,
then $S$ is \begriff{impermeable}.
Note that this property does not change under an isotopy of $S$ mod
$\t^1$.
Impermeable surfaces are closely related to (almost) $2$--normal
surfaces, by the following lemma that is implicit in the
literature~\cite{matv}, \cite{king1}. For completeness, we resume the
proof.

\begin{lemma}\label{lem:impermeable}
  Any impermeable surface $S\subset M$ with strict upper and lower
  compressing discs $D_1$ and $D_2$ is related to an almost $2$--normal
  surface with exactly one octagon by an isotopy mod $\t^1$ with support
  in $U(D_1\cup D_2)$ .
\end{lemma}
\begin{proof}
  By isotopy of $S$ mod $\t^1$ with support in $U(D_1)$, we pull the base
  of the strict upper compressing disc $D_1$ into a single tetrahedron,
  close to the string of $D_1$, changing $S$ into a surface $S_1$. Since
  $S_1$ is impermeable, no lower compressing disc for $S_1$ is contained in
  a $2$--simplex of $\t$. Hence any return of $S_1$ gives rise to an
  \emph{upper} compressing disc.
  Similarly, we pull the base of $D_2$ into a single tetrahedron by
  isotopy mod $\t^1$ with support in $U(D_2)$, changing $S$ into a
  surface $S_2$ so that any return of $S_2$ gives rise to a lower
  compressing disc.
  
  We transform $S_1$ into $S_2$ by isotopy mod $\t^1$ with support in
  $U(D_1\cup D_2)$. No surface occurring in the transformation has
  both upper and lower compressing discs in $\t^2$, since this would
  give rise to a pair of nested or independent compressing discs.
  Thus, in the course of the isotopy, we obtain a surface $\tilde S$ in
  general position to $\t$ that has neither upper nor lower
  compressing discs in $\t^2$. Hence $\tilde S$ is almost
  $k$--normal for some natural number $k$.
  It is easy to verify (see~\cite{matv} for details) that an almost
  $k$--normal surface has a pair of nested or independent compressing
  discs, if it is not almost $2$--normal or has more than one octagon.
  
  There remains to show that $\tilde S$ is not almost $1$--normal --- for
  details see Lemma~21 in~\cite{king1}. Assume that $\tilde S$ is
  almost $1$--normal. Let $D_u$ be a strict upper compressing disc for
  $\tilde S$.
  There is a subdisc $D\subset D_u$ such that $\d D\cap \tilde S$ is a
  single arc, and the closure of $\d D\setminus \tilde S$ is an arc
  $\gamma\subset \t^2$ that connects two different normal arcs of
  $\tilde S$ and is the only arc in $D\cap\t^2$ with that property.
  By pulling $\tilde S\cap \d D$ along $D$ into a single tetrahedron,
  we obtain a surface $\tilde F=\d (B^-(\tilde S)\cup U(D\setminus
  U(\gamma)))$. By the choice of $D$, $\tilde F$ has no returns, hence
  it is almost $1$--normal and $\tilde F^\times$ is isotopic mod $\t^2$
  to $\tilde S^\times$.
  The union of $D\cap U(\gamma)$ with a stripe in $\t^2$ yields a not
  necessarily strict upper compressing disc $\tilde D_u$ for $\tilde F$
  contained in a single tetrahedron $t$ whose string $\alpha$ is
  contained in $B^+(\tilde F)$. 
  We can assume that $\tilde D_u$ is contained in a single connected
  component $K$ of $t\cap B^+(\tilde F^\times)$ and meets $\d K$ only
  in its string $\alpha$. Since $\tilde F$ is almost $1$--normal, the
  two points of $\d \alpha$ belong to \emph{different} components of
  $\d t\cap \tilde F$.
  
  Since $\tilde F$ is impermeable, it has a strict lower compressing
  disc $D_l$.
  Let $\delta$ be a component of $D_l\cap \d K$. Since any component
  of $\tilde F\cap \d t$ contains at most one point of $\tilde D_u$,
  there is a subdisc $C'\subset (\d K\cap \tilde F)\setminus \t^1$
  with $\d C'\subset \tilde F\cup \delta$ and $\delta\subset \d C'$,
  such that $\d C'$ is disjoint from $D_u$. Hence, by cut-and-paste
  operations and since no lower compressing disc is completely
  contained in $K$, we can assume that $D_l$ and $\tilde D_u$ are
  independent, $D_l\cap \tilde D_u\subset \alpha$. This is impossible,
  as $\tilde F$ is impermeable. This finally proves that $\tilde S$ is
  not almost $1$--normal and has thus exactly one octagon.
\end{proof}

%***********************************************************************

\section{Upper and lower reductions}
\label{sec:reductions}

We recall here the reduction technique that we have introduced
in~\cite{king1} and refined in~\cite{king2}.
If a surface $S\subset M$ has a strict compressing disc then one can
pull $S$ along it, as in the next definition, in order to decrease the
weight $\|S\|=\#(S\cap \t^1)$.
Under certain conditions on $S$, as stated below, one eventually comes
to an almost $1$--normal surface by repeating this process. This fact is
the basis for our main applications of the reduction technique.

\begin{definition}
  Let $S\subset M$ be a splitting surface that is in general position
  to $\t^2$.  Let $D\subset M$ be a strict upper (resp.\ lower)
  compressing disc for $S$.
  The surface $\d\big(B^-(S)\cup U(D)\big)$ (resp.\ $\d\big(B^+(S)\cup
  U(D)\big)$) obtained by pulling $S$ along $D$ is an {elementary
    reduction} of $S$ along $D$.
  An upper (resp.\ lower) \begriff{reduction} of $S$ is a surface
  $S'\subset M$ obtained from $S$ by successive elementary
  reductions along strict upper (resp.\ strict lower) compressing
  discs. If these discs are contained in $\t^2$ then $S'$ is an
  upper (resp.\ a lower) $\t^2$--reduction of $S$.
\end{definition}

We do not allow mixing of elementary reductions along upper and
lower compressing discs in the transition from $S$ to $S'$.  We have
$\|S'\|\le \|S\|$ with equality if and only if $S=S'$.
The following lemma is the key tool in applications of the
reduction technique.

\begin{lemma}[Corollary~1 in~\cite{king2}, p.~57]\label{lem:t2red}
  Let $N\subset M$ be a sub--$3$--manifold such that $\d N$ is
  $k$--normal for some $k\in\mathbb N$. Let $S\subset N$ be a closed
  connected splitting surface in general position to $\t$.
  Assume that $S$ has only strict upper essential discs, and has no
  lower compressing discs contained in a single tetrahedron. Then $S$
  has an almost $1$--normal upper $\t^2$-reduction in $N$ that has
  only strict upper essential discs.\qed
\end{lemma}

The idea of the proof of the preceding lemma is as follows. Since $\d
N$ is $k$--normal, it has no returns and no circles in $\t^2\setminus
\t^1$.  Hence if a surface in $N$ has a strict compressing disc in
$\t^2$ then this disc is contained in $N$. Thus any upper
$\t^2$--reduction $S'$ of $S$ is contained in $N$.
One proves by induction on the number of elementary reductions that
$S'$ has only strict upper essential discs and has no lower compressing discs
contained in a tetrahedron.
Hence any return of $S'$ gives rise to an upper compressing disc $D\subset
\t^2$ for $S'$ with string in $B^+(S')$. If there is a circle $\gamma$ in $D\cap S'$ then
$\gamma$ bounds discs in $B^+(S')$ (e.g., a strict upper essential
disc), in \emph{both} adjacent tetrahedra. It follows that $\gamma$
and $\d D\cap S'$ belong to distinct components of $S'$. This
contradicts the hypothesis that $S$ is connected. Hence $D$ is strict
and gives rise to another elementary $\t^2$--reduction of $S'$. Thus we can
assume that $S'$ has no returns at all, i.e., $S'$ is almost
$k$--normal for some $k\in \{1,2,\ldots\}$.  Since $S'$ has no lower
compressing discs contained in a tetrahedron, $S'$ is almost
$1$--normal.

\begin{theorem}[Implicit in~\cite{jacorubinstein}, Proposition~3.3]\label{thm:gibtN}
  Let $M\not\homeo S^3$ be a closed orientable irreducible
  $3$--manifold with a triangulation $\t$.
  There is a $1$--normal sphere $F_0\subset M$ such that $B^-(F_0)$ is a
  $3$--ball and any $1$--normal sphere in $B^+(F_0)$ is isotopic mod
  $\t^2$ to $F_0$.  
\end{theorem}
\begin{proof}
  By Theorem~\ref{thm:kneserhaken}, there is a disjoint union
  $\Sigma\subset M$ of finitely many $1$--normal spheres $S_1,\dots,
  S_m$ such that any $1$--normal sphere in $M\setminus \Sigma$ is
  isotopic mod $\t^2$ to a connected component of $\Sigma$.
  The link $\d U(x)$ of a vertex $x$ of $\t$ is a $1$--normal sphere
  that can be removed from any $1$--normal surface by isotopy mod
  $\t^2$. Hence, $\Sigma$ contains a copy of the link of any vertex.
  
  Since $M$ is irreducible and not homeomorphic to $S^3$, there is a
  unique open $3$--ball $B_i\subset M$ with $\d B_i=S_i$, for
  $i=1,\dots, m$. Let $N=M\setminus \bigcup_{i=1}^m B_i$.
  For $i=1,\dots, m$, if $S_i$ is the link $\d U(x)$ of a vertex $x$
  then $x\in B_i$. Hence by definition, $N$ contains no vertex of
  $\t$. Therefore $B^+(\d N)$ does not depend on the choice of the
  vertex $x_0$ in the definition of $B^+(\cdot)$, and $N=B^+(\d N)$.
  
  We now show that $\d N$ is connected, which proves the theorem with
  $F_0=\d N$.
  Let us assume that $\d N$ is not connected, say $\d N=
  S_1\cup\dots\cup S_n$ with $n>1$. Since $N\cap \t^2$ is connected,
  there is a system $\Gamma\subset (N\cap\t^2)\setminus \t^1$ of $n-1$
  disjoint simple arcs with $\d\Gamma\subset \d N$, such that $S=(\d
  N)^\Gamma$ is a sphere. Note that $S$ is not almost $1$--normal, since
  $\Gamma$ is not in general position to $\t$.
  Since $\Gamma\subset \t^2\cap B^+(\d N)$, $S$ has only strict upper
  essential discs, and $S$ has no lower compressing discs contained in
  a single tetrahedron. Hence by Lemma~\ref{lem:t2red}, $S$ has an
  almost $1$--normal upper $\t^2$--reduction $S'\subset N$ that has
  only strict upper essential discs.
  Any connected component of $(S')^\times$ is a $1$--normal sphere in
  $N$. Hence, by the choice of $\Sigma$, $(S')^\times$ is formed by
  copies of boundary components of $N$.
  Since $S'$ is connected and has only strict upper essential discs,
  $(S')^\times$ contains \emph{at most one copy} of each connected
  component of $\d N$.  Since $S'$ does not separate two components of
  $\d N$, we have either $(S')^\times = \d N$ or
  $(S')^\times = \emptyset$.
  But since $S$ is not almost $1$--normal, it follows $S\not=S'$.
  Therefore $\|S'\| < \|S\|$, thus $(S')^\times = \emptyset$ and
  $\|S'\|=0$. By consequence, $S'$ bounds a ball $B'\subset
  M\setminus \t^1$, hence $B^+(S')=B'$. But $B^-(S')\simeq B^-(S)$ is
  a ball as well, since it is a $\d$--connected sum of the
  $3$--balls forming $M\setminus N$. This contradicts the hypothesis
  $M\not\homeo S^3$, which finally proves that $\d N$ is in fact
  connected. The theorem follows with $F_0=\d N$.
\end{proof}

We remark that, since $M\not\homeo S^3$, there is no $2$--normal
sphere in $B^+(F_0)$ with exactly one octagon (see~\cite{thompson},
\cite{matv}). A proof using the reduction technique is implicit
in~\cite{king2}, Chapter~4.
A construction algorithm for the system $\Sigma$ occurring in the
preceding proof and an estimate for $\|\Sigma\|$ can be found
in~\cite{king1}. It was applied to a study of bridge numbers of links
in~\cite{king1} and to convex $4$--polytopes in~\cite{king3}.

%**********************************************************************

\section{Completions of surfaces}
\label{sec:completion}

In this section, we introduce a bunch of notions that naturally occur
when considering strongly irreducible Heegaard surfaces.
Let $S\subset M$ be a splitting surface. A simple closed curve on $S$
contained in the boundary of a tetrahedron is a \begriff{short upper}
(resp. lower) \begriff{meridian} of $S$ if it does not bound a disc in
$S$ and does bound an embedded disc in $B^+(S)$ (resp.\ in $B^-(S)$).

Let $c \subset S$ be any simple closed curve contained in the boundary
of a tetrahedron. The curve bounds an embedded disc in the
tetrahedron. By Scharlemann's ``no nesting'' lemma (Lemma~2.2
in~\cite{scharlemann}), if $S$ is a strongly irreducible Heegaard
surface then $c$ bounds an embedded disc in $B^+(S)$ or in $B^-(S)$.
Hence, if $c$ does not bound an embedded disc in $B^-(S)$, then $c$ is a
short upper meridian of $S$. 
If $S$ is in general position to $\t$ then any pair of short meridians
can be made disjoint, by pushing one of them into the interior of a
tetrahedron. Thus, a strongly irreducible Heegaard surface in general
position to $\t$ either has no short upper meridians or no short lower
meridians.

This gives rise to the following definition. Let $S\subset M$ be a
splitting surface. If any simple closed curve on $S$ contained in the
boundary of a tetrahedron bounds an embedded disc in $B^-(S)$ (resp.
in $B^+(S)$) then $S$ is \begriff{upper trivial} (resp. lower
trivial).  If $S^\times $ is defined then we call $S$ an
\begriff{upper completion} (resp.  lower completion) of $S^\times$.
Evidently, if $S$ is upper trivial and $\Gamma\subset B^+(S)$ is a graph
in general position to $\t^2$ then $S^\Gamma$ is upper trivial as well.
By the preceding paragraph, a strongly irreducible Heegaard surface in
general position to $\t$ is upper (resp. lower) trivial if and only if
it has no short upper (resp. lower) meridian.

The following two lemmas are formulated for an upper trivial surface.
Of course, these lemmas remain true if we exchange ``upper'' with
``lower'' and $B^+$ with $B^-$.

\begin{lemma}\label{lem:ucompletionhasmeridian}
  Let $M$ be irreducible and let $S\subset M$ be a connected upper
  trivial surface. Then either $S$ has a short lower meridian or $S$ is
  isotopic to a connected component of $S^\times$ or $S$ is contained in
  a $3$--ball.
\end{lemma}
\begin{proof}
  Let us assume that $S$ has no short lower meridian. Since $S$ is
  upper trivial, any simple closed curve $c\subset S$ contained in the
  boundary of a tetrahedron $t$ bounds a disc in $B^-(S)$. By
  assumption, $c$ actually bounds a disc $D\subset S$, and it also
  bounds a disc $D'\subset t$.
  By taking subdiscs, we assume that $D\cap D'=\d D=\d D'$. Since $M$
  is irreducible, the sphere $D\cup D'$ bounds a $3$--ball $B$.  We
  change $S$ by an isotopy with support in $U(B)$, replacing $D$ with
  $D'$.
  An iteration of this process yields a surface $\tilde S\subset M$
  so that $G(\tilde S)\subset G(S)$, and any simple closed curve on
  $\tilde S$ contained in the boundary of a tetrahedron $t$ bounds a
  disc in $\tilde S\cap t$. 
  
  If $G(\tilde S)\not= \emptyset$ then $\tilde S$ is isotopic mod
  $\t^2$ to a connected component of $S^\times$.  Otherwise, $\tilde
  S$ is contained in a single tetrahedron ($\tilde S\cap
  \t^2=\emptyset$) or in the union of two tetrahedra (in this case,
  $\tilde S$ is a sphere). Hence, $S\simeq \tilde S$ is contained in a
  $3$--ball.
\end{proof}

%----------------

\begin{lemma}\label{lem:unknotb}
  Let $M$ be irreducible. Let $S\subset M$ be an almost $1$--normal
  upper trivial strongly irreducible Heegaard surface with a strict
  upper compressing disc $\tilde D$. Assume that $S$ is not isotopic
  to a connected component of $S^\times$ and $S$ is not contained in a
  $3$--ball.
  Then there is an arc $\gamma\subset (B^+(S^\times)\cap\t^2)\setminus
  \t^1$ such that $S$ is related to an upper completion of
  $(S^\times)^\gamma$ by isotopy mod $\t^1$ with support in $U(\tilde
  D)$.
\end{lemma}

\begin{proof}
  By taking a parallel copy of a subdisc of $\tilde D$, there is a
  compact embedded disc $D\subset U(\tilde D)\setminus \t^1$ in
  general position to $S$ such that the closure of $\d D\setminus S$
  is an arc $\gamma\subset B^+(S)\cap \t^2$ connecting two different
  normal arcs of $S\cap \t^2$, and no other arc in $(D\cap
  \t^2)\setminus \gamma$ connects two different normal arcs of $S$.
  When we pull $S$ across $D$ by an isotopy mod $\t^1$ with support in
  $U(\tilde D)$, we obtain the surface $S_D=\d (B^-(S)\cup U(D))$. We
  will show that $S_D$ is an upper completion of $(S^\times)^\gamma$.
  
  Consider a subdisc $D'\subset D$ such that $\d D'$ is formed by an
  arc $\alpha\subset \d t_1$ and an arc $\beta\subset S\cap t_1$, for
  some tetrahedron $t_1$. Such a subdisc exists by an innermost arc
  argument. Let $t_2$ be the other tetrahedron whose boundary contains
  $\alpha$.
  When we pull $S$ across $D'$ (as part of the transition from $S$
  to $S_D$), we obtain $S'= \d (B^-(S)\cup U(D'))$.  If
  $\alpha\not=\gamma$ then $\alpha$ does not connect different normal
  arcs of $S$, thus $(S')^\times$ is isotopic mod $\t^2$ to
  $S^\times$. If $\alpha=\gamma$ then $S'=S_D$.
  
  To prove that $S'$ is upper trivial, we have to consider the parts
  of $S'\cap \t^2$ that do not belong to $S\cap \t^2$. If $D'\cap
  \t^2$ contains simple closed curves (in the interior of $D'$) then
  these give rise to simple closed curves of $S'\cap \t^2$ bounding
  discs in $S'$.
  There remains to consider curves adjacent to $\alpha$, in the following two cases.
  
  1. We assume that $\d \alpha$ is contained in a single component of
  $(S\cap\t^2)\setminus \t^1$. Let $\alpha'\subset
  (S\cap\t^2)\setminus \t^1$ be the arc that connects the two
  endpoints of $\alpha$, and let $c_1,c_2$ be the components of $S\cap
  \d t_1,S\cap \d t_2$ that contain $\d \alpha$.
  By hypothesis and by Lemma~\ref{lem:ucompletionhasmeridian}, $S$ has
  a short lower meridian. The curves $\alpha'\cup \beta$,
  $(c_1\setminus \alpha')\cup \beta$ and $c_2$ bound discs in $M$.
  They are contained in a single tetrahedron, hence are disjoint from
  the short lower meridian.  Thus, since $S$ is strongly irreducible
  and by the ``no nesting'' lemma, they bound discs in $B^-(S)$.
  Taking the union of discs in $B^-(S)$ bounded by $c_2$ and
  $\alpha'\cup \beta$, we see that $(c_2\setminus \alpha')\cup
  \beta$ bounds a disc in $B^-(S)$ as well.
  Transforming $S$ into $S'$, $\alpha'\cup \beta$ becomes a new
  connected component of $(S'\cap\t^2)\setminus \t^1$, $(c_1\setminus
  \alpha')\cup \beta$ replaces $c_1$, and $(c_2\setminus \alpha')\cup
  \beta$ replaces $c_2$ --- compare the left part of
  Figure~\ref{fig:pull}.  Since all these curves bound discs in
  $B^-(S')$, $S'$ is upper trivial.
  \begin{figure}[htbp]
    \centering
    {\psfrag{c1'}{$c_1'$}\psfrag{c2'}{$c_2'$}\psfrag{c1}{$c_1$}\psfrag{c2}{$c_2$}\psfrag{a}{$\alpha$}\psfrag{a'}{$\alpha'$}\psfrag{D}{$D'$}\psfrag{b}{$\beta$}\psfrag{S}{$S$}\psfrag{S'}{$S'$} \epsfig{file=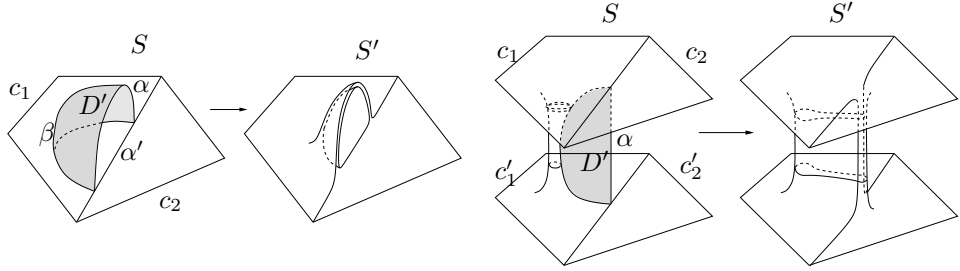,width=\textwidth}}
    \caption{The surface remains upper trivial. To keep the figure simple, the $2$--simplex containing $\alpha$ is not drawn.}
    \label{fig:pull}
  \end{figure}
  
  2. We assume that $\alpha$ connects two different components of
  $(S\cap\t^2)\setminus \t^1$. Since $S$ is almost $1$--normal,
  $\alpha$ connects two different components $c_i,c'_i$ of
  $S\cap \d t_i$, for $i=1,2$ --- see the right part of
  Figure~\ref{fig:pull}. Since $S$ is upper trivial, $c_i,c'_i$ bound
  discs in $B^-(S)$. Transforming $S$ into $S'$, these discs become
  connected along a stripe in $U(\alpha)\subset B^-(S')$, and
  therefore any connected component of $S'\cap \d t_i$ bounds a disc
  in $B^-(S')$.  Thus $S'$ is upper trivial.
  
  We replace $S$ with $S'$ and repeat until $\alpha=\gamma$. In the
  final step of this iteration, we have $S'=S_D$, hence $S_D$ is upper
  trivial. Furthermore, $G(S_D)$ is isotopic mod $\t^2$ to
  $G\left((S^\times)^\gamma\right)$. Since $S^\times$ is $1$--normal,
  $(S^\times)^\gamma\setminus \t^2$ is a disjoint union of discs.
  Hence $S_D^\times$ is isotopic mod $\t^2$ to $(S^\times)^\gamma$.
\end{proof}

Scharlemann's ``no nesting'' lemma is fundamental for our applications
of upper and lower completions.  Scharlemann's short proof (Lemma~2.2
in~\cite{scharlemann}) refers to a lemma stated
in~\cite{cassongordon}, but the proof there is also based on citation.
To avoid nested citations, we include here an elementary proof of the
``no nesting'' lemma.

\begin{lemma}[Lemma~2.2 in~\cite{scharlemann}]
  Let $H\subset M$ be a strongly irreducible Heegaard surface, and let
  $\alpha\subset H$ be a simple closed curve bounding an embedded disc
  in $M$. Then $\alpha$ bounds an embedded disc in $B^+(H)$ or
  $B^-(H)$.
\end{lemma}

\begin{proof}
  Let $D\subset M$ be an embedded disc with $\d D=\alpha$
  such that $D$ is in general position to $H$ and $\#(D\cap H)$ is
  minimal. We can assume that $D\cap H\not= \d D$, since otherwise there
  is nothing to prove. In particular, $\alpha$ is essential. We can
  assume by induction that if a simple closed curve
  $\tilde\alpha\subset H$ bounds a disc $\tilde D\subset M$
  with $\#(\tilde D\cap H) < \#(D\cap H)$ then $\tilde\alpha$ bounds
  a disc in $B^+(H)$ or $B^-(H)$.
  
  Let $\delta$ be a connected component of $(D\cap H)\setminus \d
  D$.  It bounds a subdisc $\tilde D\subset D$ with $\#(\tilde D\cap
  H) < \#(D\cap H)$. Hence by induction, $\delta$ bounds a disc in one
  handlebody of the decomposition, say, in $B^-(H)$. Cut-and-paste
  arguments and minimality of $\#(D\cap H)$ imply that 
  $\delta$ is essential on $H$ and bounds a  component of
  $D\cap B^-(H)$.
  Since $H$ is strongly irreducible, the other components of $(D\cap
  H)\setminus \d D$ are not meridians of $B^+(H)$. Hence when we apply
  the same argument, it follows that $D\cap B^-(H)$ is a non-empty
  disjoint union of meridional discs, and $P=D\cap B^+(H)$ is a
  non-empty planar surface since $D\cap H\not= \d D$ by assumption.
  
  Let $C\subset B^+(H)$ be a meridional disc transversely intersecting
  $P$ such that $\#(C\cap P)$ is minimal. Since $H$ is strongly
  irreducible and $\d P$ contains a meridian of $B^-(H)$, we have
  $C\cap P\not=\emptyset$.
  Let $\gamma\subset C\cap P$ be an innermost arc, cutting off a disc
  $C'\subset C$ with $\d C'\cap \d P =\d \gamma$. By cut-and-paste
  arguments, there is no circle in
  $C'\cap P$, hence $C'\cap P =\gamma$. Let $\beta =\d C'\cap H$ and
  $\Delta=\d P\setminus \d D$. We consider three cases.
  
  1. If $\gamma$ connects two different components of $\Delta$ then
  one can decrease $\#(D\cap H)$ by isotopy of $D$ along $C'$, in
  contradiction to the choice of $D$.
  
  2.  Assume that $\d \gamma $ is contained in a single component
  $\delta\subset \d P$. If $\gamma$ cuts off a disc $D'$ in $P$ then a
  copy of $C'\cup D'$ is a disc in $B^+(H)$ that is disjoint from $\d
  P$. Since $\d P $ contains a meridian of $B^-(H)$ and $H$ is
  strongly irreducible, $\d (C'\cup D')$ bounds a disc in $H$.  Thus
  we can decrease $\#(C\cap P)$, in contradiction to the choice of
  $C$. Hence, $\gamma$ does not cut off a disc in $P$.
  
  Let $\delta_1,\delta_2$ be the two connected components of
  $\delta\setminus\d\gamma$. If $\delta_1\cup \gamma$ (resp.
  $\delta_2\cup \gamma$) bounds a disc $D''\subset D$ then
  $D''\not\subset P$ by the preceding paragraph. Therefore $\#((C'\cup
  D'')\cap H) < \#(D\cap H)$, and hence $\delta_1\cup \beta$ (resp.
  $\delta_2\cup \beta$) bounds an embedded disc in $B^+(H)$ or
  $B^-(H)$ by induction. It actually bounds a disc in $B^-(H)$, since
  $\d P$ contains a meridian of $B^-(H)$, since a copy of
  $\delta_i\cup \gamma$ is disjoint from $\d P$, and since $H$ is
  strongly irreducible.
  
  If $\delta = \d D$ then both $\delta_1\cup \gamma$ and $\delta_2\cup
  \gamma$ bound discs in $D$. Hence by the preceding paragraph, $\d D$
  bounds a disc in $B^-(H)$, in contradiction to our assumption. If
  $\delta\subset \Delta$ then either $\delta_1\cup \gamma$ or
  $\delta_2\cup \gamma$ bounds a disc $D'\subset D$. Hence it bounds a
  disc $D''\subset B^-(H)$. Replacing $D'$ with $D''\cup C'$, we can
  decrease $\#(D\cap H)$, since $D'\not\subset P$. We obtain a
  contradiction to the minimality of $\#(D\cap H)$.
  
  3. We are left with the case that $\gamma$ connects $\d D$ with a
  component $\delta$ of $\Delta$. Let $D_\delta$ be the disc in $D\cap
  B^-(H)$ bounded by $\delta$. Let $\tau\subset H$ be a simple closed
  curve obtained by connecting $\delta\setminus U(\d \beta)$ and $\d
  D\setminus U(\d \beta)$ with two copies of $\beta$. By taking the
  union of the disc $D\setminus (U(\gamma)\cup D_\delta)$ with two
  copies of $C'$, $\tau $ bounds an embedded disc $\tilde D\subset M$
  with $\#(\tilde D\cap H) < \#(D\cap H)$. Hence, using the same
  arguments as in Case 2, $\tau$ bounds a disc $D'$ in $B^-(H)$.
  Taking the union of $D'$ with a stripe along $\beta$ and with
  $D_\delta$, we obtain a disc in $B^-(H)$ bounded by $\d D$.
  This contradicts our assumption and proves the lemma.
\end{proof}

%*************************************************
%*************************************************

\section{Proof of Theorem~\ref{thm:heegaard-an}}
\label{sec:proof}

Let $M$ be a closed orientable irreducible $3$--manifold with a
triangulation $\t$, and let $H\subset M$ be a strongly irreducible
Heegaard surface of $M$.
By Waldhausen's theorem~\cite{waldhausen} on Heegaard surfaces of
$S^3$, $H$ is an embedded $2$--sphere if $M\homeo S^3$.  In this case,
we pick a vertex $x\in \t^0$, and connect two copies of the
$1$--normal sphere $\d U(x)$ along an unknotted arc. We obtain an
almost $1$--normal sphere isotopic to $H$ and satisfying
Theorem~\ref{thm:heegaard-an}.
There remains to consider the case $M\not\homeo S^3$.  In particular, $H$
is not contained in a $3$--ball. Let $F_0\subset M$ be as in
Theorem~\ref{thm:gibtN}, and let $N=B^+(F_0)$.

By isotopy of $H$, we can assume that the ball $M\setminus N$ is
contained in a single connected component $\mathcal H$ of $M\setminus
H$. Since $\mathcal H$ is a handlebody, there is a graph
$\Gamma_-\subset \mathcal H\cap N$ in general position to $\t$ with $\d
\Gamma_-\subset \d N$, such that $\mathcal H$ collapses onto $(M\setminus
N)\cup \Gamma_-$. Hence $S_-=(\d N)^{\Gamma_-}$ is an almost
$1$--normal upper trivial surface isotopic to $H$.
Similarly, there is a graph $\Gamma_+\subset M\setminus \mathcal H$ in
general position to $\t$ such that $M\setminus \mathcal H$ collapses
onto $\Gamma_+$. Hence $S_+=\d U(\Gamma_+)$ is an almost
$1$--normal lower trivial surface isotopic to $H$.
We have $N\supset B^+(S_-)\supset B^+(S_+)= U(\Gamma_+)$.
  
Thus there are $1$--normal surfaces $F_-,F_+\subset N$, an upper
completion $S_-\simeq H$ of $F_-$ and a lower completion $S_+\simeq H$
of $F_+$ with $B^+(S_+)\subset B^+(S_-)$. Let $S\simeq H$ be a lower
(resp. upper) trivial almost $1$--normal surface that is nested
between $S_-$ and $S_+$, i.e., $B^+(S_+)\subset B^+(S)\subset
B^+(S_-)$.
Assume that $\|S\|<\|S_-\|$ (resp. $\|S\|<\|S_+\|$) or there is a
connected component of $S^\times$ that can not be isotoped mod $\t^2$
into $B^-(F_-)$ (resp.  of $B^+(F_+)$). Then we replace $F_-$ by
$S^\times$ (resp. $F_+$ by $S^\times$), and iterate.
By Kneser's finiteness result (Theorem~\ref{thm:kneserhaken}), the
iteration stops after a finite number of steps. We choose $F_-,F_+$ so
that a further iteration is impossible.

If some component $F$ of $F_-$ or $F_+$ is isotopic to $H$ then we
connect $F$ with $\d N$ along an unknotted arc, and obtain the first
alternative of Theorem~\ref{thm:heegaard-an}. By now we assume that no
connected component of $F_-$ or $F_+$ is isotopic to $H$. Hence by
Lemma~\ref{lem:ucompletionhasmeridian}, $S_-$ has a short lower
meridian and $S_+$ has a short upper meridian.

\smallskip
\paragraph*{\textbf{Case 1.}}
We assume that $S_-$ has a strict upper compressing disc in $B^-(S_+)$, for
some choice of $S_-,S_+$.
According to Lemma~\ref{lem:unknotb}, $S_-$ is isotopic mod $\t^1$ to an
upper completion of $F_-^\gamma$ contained in $B^+(S_-)\cap B^-(S_+)$,
for some arc $\gamma\subset (B^+(F_-)\cap \t^2)\setminus \t^1$.
The surface $F_-^\gamma$ satisfies the hypothesis of
Lemma~\ref{lem:t2red}, hence $F_-^\gamma$ has an almost $1$--normal
upper $\t^2$-reduction $S'\subset B^+(F_-)\cap B^-(F_+)$ with $\|S'\|<
\|F_-\|$.
We use the freedom in the choice of $S_-,S_+$ and isotope it so that
$S_-\cup S_+$ has empty intersection with the strict compressing discs
involved in the elementary reductions transforming $F_-^\gamma$ into
$S$. In particular, the upper $\t^2$--reduction $S'$ of $F_-^\gamma$
gives rise to an almost $1$--normal upper $\t^2$--reduction
$S'_-\subset B^+(S_-)\cap B^-(S_+)$ of $S_-$.

If $F_-^\gamma\simeq H$ then we get the first alternative of
Theorem~\ref{thm:heegaard-an}, by pushing $\gamma$ into the interior of
a tetrahedron. Otherwise, by Lemma~\ref{lem:ucompletionhasmeridian},
$S_-$ has a small lower meridian. We take a copy $c\subset S_-\setminus
\t^2$ of this meridian.  It is not affected by elementary reductions
along compressing discs in $\t^2$. Thus $S'_-$ has no short upper
meridian, since it would be disjoint from the meridian $c$ of $B^-(S'_-)$.
Hence $S'_-\simeq H$ is upper trivial.  This is a contradiction to the
choice of $F_-,F_+$, as $\|S'_-\|< \|S_-\| = \|F_-\|$.

In conclusion, if $S_-$ has a strict upper compressing disc in
$B^-(S_+)$, or similarly if $S_+$ has a strict lower compressing disc
in $B^+(S_-)$, for some choice of $S_-,S_+$, then the first
alternative of Theorem~\ref{thm:heegaard-an} holds.

\smallskip
\paragraph*{\textbf{Case 2.}}
We assume that $S_-$ has no strict upper compressing disc in
$B^-(S_+)$ and $S_+$ has no strict lower compressing disc in
$B^+(S_-)$, for \emph{any} choice of $S_-,S_+$.
Let $\Sigma_g$ be the closed orientable surface whose genus $g$
coincides with the genus of $H$. By hypothesis on $F_-,F_+$, there is
an embedding $J\co \Sigma_g\times [0,1] \to M$ in general position to
$\t$ such that $J_\xi = J(\Sigma_g\times \{\xi\})\simeq H$ for all
$\xi\in [0,1]$, $J_{0}$ is an upper completion of $F_-$, and $J_{1}$
is a lower completion of $F_+$.

For any $\xi \in [0,1]$, $J_\xi$ contains no vertex ($\t^0$ is in the
complement of $N$), and it has at most one point of tangency with
$\t^1\setminus \t^0$ or $\t^2\setminus \t^1$ by general position of
$J$. A parameter $\xi\in [0,1]$ is \begriff{critical} if $J_\xi$ is
tangent to $\t^1$. There are only finitely many
critical parameters.
The \begriff{complexity} of $J$ is
$$
\kappa(J) = \sum_{\xi\in I \text{ critical}} \|J_\xi\|.$$
We choose $J$ so that $\kappa(J)$ is minimal among all embeddings
with the above properties. This is called \emph{thin position}, a
notion with many fruitful applications --- for instance in~\cite{gabai},
\cite{thompson}, \cite{matv}, \cite{scharlthompson} and~\cite{king1}.

By assumption on $F_-,F_+$, there is a short lower meridian of $J_0$
and a short upper meridian of $J_1$. Hence there is a non-critical
parameter $\zeta\in[0,1]$ such that $J_{\zeta}$ has neither short
upper nor short lower meridians, or has \emph{both} short upper and
short lower meridians.
But the latter case can not occur since otherwise by small isotopy we
obtain meridians of $B^-(J_{\zeta})$ and $B^+(J_{\zeta})$ that
intersect in at most one point (namely in a point of tangency of
$J_{\zeta}$ with $\t^2$), in contradiction to strong irreducibility of
$J_{\zeta}$.
Thus $J_{\zeta}$ is both upper and lower trivial, and by small isotopy we
can assume that it is in general position to $\t$.

If $J$ has no critical parameter then it easily follows that
$J_{\zeta}$ has no returns, since $J_0$ and $J_1$ have no returns,
$G(J_0)\subset B^-(J_\zeta)$, $G(J_1)\subset B^+(J_\zeta)$, and
$G(J_0)$ is isotopic mod $\t^2$ to $G(J_1)$.
Thus $J_\zeta$ is almost $k$--normal for some $k\in \mathbb N$. Since
a $k$--normal surface is determined by its intersection with $\t^1$,
$J_{\zeta}^\times$ is isotopic mod $\t^2$ to $F_-$. Hence $J_\zeta$ is
an upper completion of $F_-$ isotopic to $H$ without a short lower
meridian, in contradiction to the assumption.
Therefore, $J$ has a critical parameter.

Since we are in Case~2, $J_1=S_+$ has no strict lower compressing discs.
Hence there is a smallest critical parameter $\xi_0\in [0,1]$ with
$\|J_{\xi_0+\epsilon}\|<\|J_{\xi_0-\epsilon}\|$ for small
$\epsilon>0$.
Since $J_0$ has no strict upper compressing disc and $\xi_0$ is
minimal, we have $\|J_{0}\| < \|J_{\xi_0-\epsilon}\|$.
Thus $J_{\xi_0-\epsilon}$ has both strict upper and lower compressing
discs in $B^+(J_0)\cap B^-(J_1)$.

Assume that $J_{\xi_0-\epsilon}$ has a pair of nested or independent
upper and lower compressing discs $D_1,D_2\subset M$.
By Lemma~\ref{lem:discparallel} and since $F_-,F_+$ are $1$--normal,
we can choose $D_1,D_2\subset B^+(F_-)\cap B^-(F_+)$.
We can assume that the intersection of $D_i$ with $S_\pm$ is disjoint
from $\t^2$. Thus $D_1\cup D_2$ does not contain a meridian of
$B^+(S_-)$ or of $B^-(S_+)$, by the existence of short meridians of
$S_\pm$ and strong irreducibility.
Hence an isotopy of $J$ along $D_1 \cup D_2$ yields an embedding
$J'\co \Sigma_g\times I\to M$ so that $J'_0$ is an upper completion of
$F_-$ and $J'_1$ is a lower completion of $F_+$. Moreover, the number
of critical parameters does not increase and the weight of some
critical levels \emph{does} decrease, which implies
$\kappa(J')<\kappa(J)$ --- see~\cite{matv} or \cite{thompson} for
details.  This contradicts the choice of $J$ and disproves the
existence of $D_1,D_2$.

By the preceding paragraphs, $J_{\xi_0-\epsilon}\subset B^+(J_0)\cap
B^-(J_1)$ is impermeable. It has strict upper and lower compressing
discs in $B^+(J_0)\cap B^-(J_1)$ yield by critical points of $J$.
Thus by Lemma~\ref{lem:impermeable}, $J_{\xi_0-\epsilon}$ is isotopic
mod $\t^1$ to an impermeable almost $2$--normal surface $S\subset
B^+(J_0)\cap B^-(J_1)$ that has exactly one octagon.

We assume that $S$ has a short lower meridian; the opposite case
of a short upper meridian is of course symmetric.
Let $F$ be the component of $S^\times$ containing the octagon. The
octagon yields a strict upper compressing disc for $F$ contained in a
tetrahedron. Let $F'$ be the result of an elementary reduction of $F$
along that strict upper compressing disc.
$F'$ satisfies the hypothesis of Lemma~\ref{lem:t2red}. Thus $F'$ has
an almost $1$--normal upper $\t^2$--reduction $F''\subset B^+(S^\times)\cap
B^-(F_+)$.
As in Case~1, by an appropriate choice of $S_-$ and $S_+$, we obtain
an almost $1$--normal upper reduction $S'\subset B^+(S_-)\cap
B^-(S_+)$ of $S$ that has a short lower meridian. Hence $S'$ is upper
trivial.
The $1$--normal surface $(S')^\times$ is separated from $F_-$ by $F$,
which has an octagon. Hence some connected component of $(S')^\times$
can not be isotoped mod $\t^2$ into $B^-(F_-)$. This is a
contradiction to the choice of $F_-$.

Thus $S$ has neither short upper nor short lower meridians, i.e., $S$
is both upper and lower trivial. Therefore $H$ is isotopic to a
connected component $F_1$ of $S^\times$, by
Lemma~\ref{lem:ucompletionhasmeridian}. If $F_1$ contains the octagon
then we get the second alternative of Theorem~\ref{thm:heegaard-an}.
Otherwise, we join $F_1$ with $\d N$ along an unknotted arc and get
the first alternative of Theorem~\ref{thm:heegaard-an}.
This finishes the proof of Theorem~\ref{thm:heegaard-an}.


\begin{thebibliography}{[14]}
\bibitem{cassongordon}
\newblock Casson, A.\ and Gordon, C.\ McA.:
\newblock Reducing Heegaard splittings.
\newblock {\em Topology Appl.} \textbf{27,} 275--283 (1987).

\bibitem{gabai}
\newblock Gabai, D.:
\newblock Foliations and the topology of 3-manifolds III.
\newblock {\em J. Differ. Geom.} \textbf{26,} 445--503 (1987).

%\bibitem{haken1}
%\newblock Haken, W.:
%\newblock {\"U}ber das Hom{\"o}omorphieproblem der
%3-Mannigfaltigkeiten I. 
%\newblock {\em Math. Z.} \textbf{80,} 89--120 (1962).

\bibitem{haken2}
\newblock Haken, W.:
\newblock \textit{Some results on surfaces in
  3-manifolds}. In: Studies in Modern Topology (MAA Stud.
Math. 5, Prentice Hall 1968) 39--98.
\newblock Editor Hilton,~P.

%\bibitem{hemion}
%\newblock Hemion, G.:
%\newblock \textit{The Classification of Knots and 3-Dimensional Spaces}.
%\newblock (Oxford University Press, 1992).


\bibitem{jacorubinstein}
\newblock Jaco, W.\ and Rubinstein, J.\ H.:
\newblock 0-efficient triangulations of 3-manifolds.
\newblock Preprint 2002.
\texttt{http://de.arXiv.org/abs/math.GT/0207158}


\bibitem{johannson}
\newblock Johannson, K.:
\newblock \textit{Topology and combinatorics of 3-manifolds.}
\newblock Lecture Notes in Mathematics 1599
\newblock (Springer-Verlag 1995).

\bibitem{king1}
\newblock King, S.:
\newblock The size of triangulations supporting a given link.
\newblock {\em Geometry \& Topology} \textbf{5,} 369--398 (2001).

\bibitem{king2}
  \newblock King, S.:
  \newblock {\em Polytopality of  triangulations}.
  \newblock PhD thesis, University of Strasbourg (2001).\\
  \texttt{http://www-irma.u-strasbg.fr/irma/publications/2001/01017.shtml.} 

\bibitem{king3}
\newblock King, S.:
\newblock How to make a triangulation of $S^3$ polytopal.
\newblock \texttt{arxiv:math.GT/0009216}. To appear in {\em Trans. Amer. Math. Soc}.

\bibitem{kneser}
\newblock Kneser, H.:
\newblock Geschlossene Fl\"achen in dreidimensionalen
Mannigfaltigkeiten.
\newblock {\em Jahresber. Deut. Math. Ver.} \textbf{38,} 248--260,
(1929). 

\bibitem{matv}
\newblock Matveev, S. V.:
\newblock {An algorithm for the recognition of 3-spheres (according to
  Thompson)}. 
\newblock {\em Mat. sb.} \textbf{186,5,} 69--84 (1995).
\newblock English translation in {\em Sb. Math.} \textbf{186,}
(1995). 

%\bibitem{matveev2}
%\newblock Matveev, S. V.:
%\newblock On the recognition problem for {H}aken 3-manifolds. 
%\newblock {\em Suppl. Rend. Circ. Mat. Palermo.} \textbf{49,} 131--148
%(1997). 


\bibitem{rubinstein}
\newblock Rubinstein, J. H.:
\newblock {\em Polyhedral minimal surfaces, {H}eegaard splittings and
  decision problems for 3-dimensional manifolds}.
\newblock Proc. Georgia Topology Conf. (Amer. Math. Soc./ Intl. Press,
1993). 

\bibitem{scharlemann}
\newblock Scharlemann, M.:
\newblock Local detection of strongly irreducible Heegaard splittings. 
\newblock {\em Topology Appl.} \textbf{90,} 135--147 (1998).

\bibitem{scharlthompson}
\newblock Scharlemann, M. and Thompson, A.:
\newblock Thin position for $3$-manifolds.
\newblock {\em Contemp. Math.} \textbf{164,} 231--238
Geometric topology (Haifa, 1992).

\bibitem{stocking}
\newblock Stocking, M.:
\newblock Almost normal surfaces in $3$-manifolds.
\newblock {\em Trans. Amer. Math. Soc.} \textbf{352,}  171--207 (2000).

\bibitem{thompson}
\newblock Thompson, A.:
\newblock Thin position and the recognition problem for $S^3$.
\newblock {\em Mathematical Research Letters} \textbf{1,} 613--630
(1994). 

\bibitem{waldhausen}
\newblock Waldhausen, F.:
\newblock Heegaard--Zerlegungen der $3$--Sph\"are.
\newblock {\em Topology} \textbf{7,}  195--203 (1968).

\end{thebibliography}
\end{document}